\documentclass[12pt]{article}

\usepackage[T1]{fontenc}
\usepackage{avant}
\usepackage[cp1250]{inputenc}
\usepackage{amsmath,amssymb,amsthm}
\usepackage{caption}
\usepackage{geometry}
\usepackage{float}

\usepackage{makeidx}
\usepackage[matrix,arrow,curve]{xy}
\usepackage{boxedminipage}
\usepackage{epic}
\usepackage{longtable}
\usepackage{ifthen}
\usepackage{tabularx}
\usepackage{parskip}
\usepackage{verbatim}
\usepackage{multirow}

\usepackage{authblk}

\usepackage{xcolor}

\newtheorem{theorem}{Theorem}
\newtheorem{proposition}{Proposition}
\newtheorem{lemma}{Lemma}

\newtheorem{remark}{Remark}

\title{On semisimple standard compact Clifford-Klein forms}

\author{Maciej Boche\'nski\\
 Faculty of Mathematics and Computer Science,\\
  University of Warmia and Mazury in Olsztyn, \\
  S\l\/oneczna 54, 10-710, Olsztyn, Poland\\
  e-mail: mabo@matman.uwm.edu.pl \\ }

\begin{document}
\maketitle

\begin{abstract}
In this paper we give the classification of standard compact Clifford-Klein forms corresponding to triples $(\mathfrak{g}, \mathfrak{h}, \mathfrak{l})$ such that $\mathfrak{g}=\mathfrak{h}+\mathfrak{l}$ and $\mathfrak{g}$ is a sum of two absolutely simple ideals. The classification is done using Onishchik's results concerning semisimple decompositions of semisimple Lie algebras. Using this classification we obtain new examples of reductive homogeneous spaces admitting non-standard compact Clifford-Klein forms.

\end{abstract}

\section{Introduction}
Let $G$ be a non-compact real linear semisimple Lie group and $H$ a reductive subgroup such that $G/H$ is of reductive type. We say that the homogeneous space $G/H$ admits a \textit{compact Clifford-Klein form} if there exists a discrete subgroup \textcolor[rgb]{0,0,1}{$\Gamma \subset G$} such that $\Gamma$ acts freely, properly and co-compactly on $G/H.$ Notice that in a case when $H$ and $G/H$ are non-compact such $\Gamma$ may not exist. The problem of finding and classifying homogeneous spaces admitting compact Clifford-Klein forms is in general difficult and is an active area of study (see for example \cite{btc, bo, kas, ko, kdef, ku, ok}).

There is basically one known way of constructing non-trivial examples of compact Clifford-Klein forms. Assume there exists a reductive connected subgroup $L\subset G$ such that $L$ acts properly and co-compactly on $G/H.$ Then one can choose a co-compact lattice $\Gamma$ of $L$ so that $\Gamma$ acts properly, co-compactly and freely on $G/H.$ In such case we say that $G/H$ admits a \textit{standard compact Clifford-Klein form}. The notion of standard compact Clifford-Klein forms was introduced by Kobayashi in \cite{ko-ma}, who also gave a list of examples of simple symmetric spaces admitting standard compact Clifford-Klein forms in \cite{kdis} and \cite{ko} (the terminology \textit{standard} was first used in \cite{kako}). There is a conjecture concerning standard Clifford-Klein forms (see Conjecture 3.3.10 in \cite{ko}) stating that if $G/H$ admits a compact Clifford-Klein form then $G/H$ also admits a standard one.

\noindent
The classification of irreducible simple symmetric spaces which admit standard compact Clifford-Klein forms was recently given in \cite{tojo}. \textcolor[rgb]{0,0,1}{In \cite{bo} it is shown that $G/H$ does not admit standard compact Clifford-Klein forms when $G$ is a simple exceptional Lie group.}

Assume that the Lie algebra $\mathfrak{g}$ of $G$ is a direct sum of two absolutely simple ideals (that is, ideals wit simple complexifications), $\mathfrak{g}_i, \ i=1,2;$ 
$$\mathfrak{g}=\mathfrak{g}_1 \oplus \mathfrak{g}_{2} .$$
Denote by $\mathfrak{h}$ the Lie algebra of a connected reductive subgroup $H\subset G.$ In this paper we list all possible triples $(\mathfrak{g}, \mathfrak{h}, \mathfrak{l})$ such that $\mathfrak{g}=\mathfrak{h}+\mathfrak{l}$ and $\mathfrak{l}$ is a Lie algebra of a connected reductive subgroup $L\subset G$ that acts properly and co-compactly on $G/H.$ This gives the list of all homogeneous spaces $G/H$ which admit standard compact Clifford Klein forms such that $G=HL$ and the lie algebra of $\mathfrak{g}$ is a sum of two absolutely simple ideals. \textcolor[rgb]{0,0,1}{Notice that in \cite{bo, tojo} triples $(G,H,L)$ corresponding to arbitrary standard compact Clifford-Klein forms were consider. On the other hand in this paper we restrict our attention to the case when $L$ acts transitively on $G/H.$ It is not known if co-compactness of the action of $L$ on $G/H$ implies that $L$ acts transitively on $G/H.$}

Notice that we may assume that $\mathfrak{g}_{1}, \mathfrak{g}_{2}$ are of non-compact type. Also we may assume that $\mathfrak{h},$ $\mathfrak{l}$ are semisimple without ideals of compact type. This is because the result of \cite{bela} which states that if a reductive homogeneous space of a semisimple Lie group admits a compact Clifford-Klein form then the center of the isotropy subgroup is compact. In other words if $\mathfrak{h}=\mathfrak{c}\oplus \mathfrak{h}_{1}$ where $\mathfrak{c}$ is an ideal of compact type then $(\mathfrak{g},\mathfrak{h},\mathfrak{l})$ induces a standard Clifford-Klein form if and only if $(\mathfrak{g},\mathfrak{h}_{1},\mathfrak{l})$ induces a standard Clifford-Klein form.

Main results of the paper are given in Theorem \ref{twg}, Table \ref{tab1} and in Proposition \ref{prop}. The proof of Theorem \ref{twg} is based on the classification of semisimple decompositions of simple and semisimple Lie algebras given in \cite{on1, on2}. The proof of Proposition \ref{prop} is carried out using the ideas and results of Kobayashi \cite{koc, kdef, ko-ma} and Kassel \cite{kas}.

Notice that it is reasonable to expect that  the class of standard Clifford-Klein forms of semisimple groups is significantly richer than the class of standard Clifford-Klein forms of simple Lie groups. However,    
in this article we present a rather surprising result that the class of standard compact Clifford-Klein forms of semisimple groups is rather small. It suggests that basically one should not expect to greatly enrich the class of compact Clifford-Klein forms beyond the double quotients $\Gamma\backslash G/H$ of simple Lie groups $G$.  

\noindent
In fact, the only non-trivial \textcolor[rgb]{0,0,1}{(i.e. neither $\mathfrak{h}$ nor $\mathfrak{l}$ contains a non-zero ideal of $\mathfrak{g}$ and $\mathfrak{g}\neq (\mathfrak{g}_{1}\cap \mathfrak{h} + \mathfrak{g}_{1}\cap \mathfrak{l})\oplus (\mathfrak{g}_{2}\cap \mathfrak{h} + \mathfrak{g}_{2}\cap \mathfrak{l})$)} semisimple compact Clifford-Klein forms such that $\mathfrak{g}=\mathfrak{g}_{1}\oplus \mathfrak{g}_{2}=\mathfrak{h}+\mathfrak{l}$ are the following

\begin{enumerate}
	\item $G\times G/\textrm{diag}(G),$ $L=G'\times G'',$ where the triple $(G,G',G'')$ is given by a triple in Table \ref{tab1}.
	\item $SO(4,4)\times SO(3,4)/\Delta (SO(3,4)),$ $L=Spin(3,4)\times SO(1,4), $ where 
	$$\Delta (SO(3,4)):=\{ (\iota_{1}(g),g) \ | \ g\in SO(3,4) \} \  \textrm{for} \ \iota_{1}:SO(3,4)\hookrightarrow SO(4,4) .$$
	\item $SO(4,4)\times SO(2,4)/\Delta (SO(2,4)),$ $L=Spin(3,4)\times SO(1,4), $ where 
	$$\Delta (SO(2,4)):=\{ (\iota_{2}(g),g) \ | \ g\in SO(2,4) \} \  \textrm{for} \ \iota_{2}:SO(2,4)\hookrightarrow SO(4,4) .$$
	\item $SO(3,4)\times SO(2,4)/\Delta (SO(2,4)),$ $L=G_{2(2)}\times SO(1,4), $ where 
	$$\Delta (SO(2,4)):=\{ (\iota_{3}(g),g) \ | \ g\in SO(2,4) \} \  \textrm{for} \ \iota_{3}:SO(2,4)\hookrightarrow SO(3,4) .$$
\end{enumerate}
\begin{remark}
In the item $(k+1)$ in the above list the embedding of $\mathfrak{l}$ into $\mathfrak{g}_{1}\oplus \mathfrak{g}_{2}$ \textcolor[rgb]{0,0,1}{is given by $\iota''\oplus j_{k},$ $k=1,2,3$ where $j_{k}:\mathfrak{g}_{2}\cap\mathfrak{l}\rightarrow \mathfrak{g}_{2}$ is s.t. 
$$\iota' = \iota_{k}\circ j_{k},$$}
and $\iota' , \iota''$ are given in Table \ref{tab1} (in the fourth row for cases 2., 3., and in the fifth row for case 4.).
\end{remark}

However, one may find in the above list interesting examples of reductive homogeneous spaces admitting non-standard Clifford-Klein forms. For the spaces listed in part 1. one may use the following result to obtain such forms
\begin{proposition}[\cite{kako}, Proposition 2.1]
Any standard compact Clifford-Klein form $\Gamma' \times \Gamma'' \setminus G\times G / \textrm{diag}(G),$ where $\Gamma'\subset G',$ $\Gamma''\subset G''$ are irreducible uniform lattices, remains a manifold after any small deformation of $\Gamma' \times \Gamma''$ inside $G\times G.$
\label{propkako}
\end{proposition}

\noindent
Using the above result we prove that
\begin{proposition}
The spaces $SO(4,4)\times SO(2,4)/\Delta (SO(2,4)),$ $SO(3,4)\times SO(2,4)/\Delta (SO(2,4))$ admit non-standard compact Clifford Klein-forms. 
\label{prop}
\end{proposition}

  \vskip6pt
   \noindent {\bf Acknowledgment.} The  author acknowledges the support of the National Science Center (grant NCN no. 2018/31/D/ST1/00083).

\section{Results and proofs}

We start with the following observation. If $\mathfrak{g}=\mathfrak{h}+\mathfrak{l}$ then (by Theorem 3.1 in \cite{on1}) we have $G=HL.$ Therefore
$$G/H \cong L/H\cap L$$
as a L-manifold so $L$ acts properly (and co-compactly) on $G/H$ if and only if $H\cap L$ is compact (a short proof of the fact that $\mathfrak{g}=\mathfrak{h}+\mathfrak{l}$ implies $G/H \cong L/H\cap L$ can be found in \cite{koinv}, Lemma 5.1). In order to find all triples $(\mathfrak{g},\mathfrak{h},\mathfrak{l})$ corresponding to standard Clifford-Klein forms we have to list all semisimple decompositions $\mathfrak{g}=\mathfrak{h}+\mathfrak{l}$ such that $\mathfrak{h}\cap \mathfrak{l}$ is of compact type. Denote by 
$$\pi_{i}:\mathfrak{g}\rightarrow \mathfrak{g}_{i}, \ \ i=1,2, $$
the projection on i-th simple factor of $\mathfrak{g}.$
\begin{theorem}
Let $\mathfrak{g}=\mathfrak{g}_{1}\oplus\mathfrak{g}_{2}$ be a semisimple Lie algebra which is a sum of two absolutely simple ideals of non-compact type. Let $\mathfrak{h},\mathfrak{l}\subset \mathfrak{g}$ be two semisimple subalgebras without non-trivial ideals of compact type, such that $\mathfrak{g}=\mathfrak{h}+\mathfrak{l} .$ Then the triple $(\mathfrak{g},\mathfrak{h},\mathfrak{l})$ gives a standard compact Clifford-Klein form if and only if it is contained in the following list
\begin{enumerate}
	\item $\mathfrak{g}_{1},\mathfrak{g}_{2}$ are arbitrary absolutely simple Lie algebras,  $\mathfrak{h}= \mathfrak{g}_{1}\oplus \{ 0 \} ,$ $\mathfrak{l}= \{ 0 \} \oplus \mathfrak{g}_{2}$ and 
	$$\pi_{1}(\mathfrak{h})=\mathfrak{g}_{1}, \ \pi_{2}(\mathfrak{h})=\{ 0 \} , \  \pi_{1}(\mathfrak{l})=\{ 0 \} , \ \pi_{2}(\mathfrak{l}) = \mathfrak{g}_{2} .$$	
	\item $\mathfrak{g}_{2}$ is an arbitrary absolutely simple Lie algebra, $\mathfrak{h}= \mathfrak{g}_{1}'\oplus \{ 0 \} ,$ $ \mathfrak{l}= \mathfrak{g}_{1}''\oplus \mathfrak{g}_{2},$ the triple $(\mathfrak{g}_{1},\mathfrak{g}_{1}',\mathfrak{g}_{1}'')$ is contained in Table \ref{tab1} and 
	$$\pi_{1}(\mathfrak{h})=\mathfrak{g}_{1}', \ \pi_{2}(\mathfrak{h})=\{ 0 \} , \  \pi_{1}(\mathfrak{l})=\mathfrak{g}_{1}'' , \ \pi_{2}(\mathfrak{l}) = \mathfrak{g}_{2} .$$
	\item $\mathfrak{h}= \mathfrak{g}_{1}'\oplus \mathfrak{g}_{2}'$, $\mathfrak{l}= \mathfrak{g}_{1}''\oplus \mathfrak{g}_{2}''$ where the triples $(\mathfrak{g}_{1},\mathfrak{g}_{1}',\mathfrak{g}_{1}''),$ $(\mathfrak{g}_{2},\mathfrak{g}_{2}',\mathfrak{g}_{2}'')$ are contained in Table \ref{tab1} and
	$$\pi_{1}(\mathfrak{h})=\mathfrak{g}_{1}', \ \pi_{2}(\mathfrak{h})=\mathfrak{g}_{2}' , \  \pi_{1}(\mathfrak{l})= \mathfrak{g}_{1}'', \ \pi_{2}(\mathfrak{l}) = \mathfrak{g}_{2}'' .$$
		\item $\mathfrak{g}_{1},\mathfrak{g}_{2}$ are arbitrary absolutely simple Lie algebras such that $\mathfrak{g}_{2}\subset\mathfrak{g}_{1}$, $\mathfrak{h}= \mathfrak{g}_{1}\oplus \{ 0 \} ,$ $ \mathfrak{l}\cong \mathfrak{g}_{2}$ and 
	$$\pi_{1}(\mathfrak{h})=\mathfrak{g}_{1}, \ \pi_{2}(\mathfrak{h})=\{ 0 \} , \  \pi_{1}(\mathfrak{l})= \mathfrak{g}_{2}, \ \pi_{2}(\mathfrak{l}) = \mathfrak{g}_{2} .$$
	\item Take a triple $(\mathfrak{g}_{2},\mathfrak{g}',\mathfrak{g}'')$ from Table \ref{tab1}. Let $\mathfrak{g}_{1}$ be an arbitrary absolutely simple Lie algebra such that $\mathfrak{g}''\subset\mathfrak{g}_{1}$. Put $\mathfrak{h}= \mathfrak{g}_{1}\oplus \mathfrak{g}' ,$ $ \mathfrak{l}\cong \mathfrak{g}''$ and 
	$$\pi_{1}(\mathfrak{h})=\mathfrak{g}_{1}, \ \pi_{2}(\mathfrak{h})=\mathfrak{g}' , \  \pi_{1}(\mathfrak{l})= \mathfrak{g}'', \ \pi_{2}(\mathfrak{l}) = \mathfrak{g}'' .$$
	\item Take a triple $(\mathfrak{g}_{1},\mathfrak{g}',\mathfrak{g}'')$ from Table \ref{tab1}. Take any simple subalgebra $\mathfrak{g}_{2}\subset \mathfrak{g}_{1}$  so that $\mathfrak{g}''\subset \mathfrak{g}_{2}\subset \mathfrak{g}_{1}.$ Put $\mathfrak{h}=\mathfrak{g}'\oplus \mathfrak{g}''$ and $\mathfrak{l}=\mathfrak{g}_{2}$ so that
	$$\pi_{1}(\mathfrak{h})=\mathfrak{g}', \ \pi_{2}(\mathfrak{h})=\mathfrak{g}'' , \ \pi_{1}(\mathfrak{l})= \mathfrak{g}_{2}, \ \pi_{2}(\mathfrak{l}) = \mathfrak{g}_{2} .$$
	All triples of that type are listed in Table \ref{tab2}.
\end{enumerate}
\label{twg}
\end{theorem}
\begin{remark}
Given a triple $(\mathfrak{g},\mathfrak{h},\mathfrak{l})$ and corresponding groups $G,H,L$ we see that if $G/H$ admits a standard compact Clifford-Klein form then $G/L$ admits a standard compact Clifford-Klein form. Therefore one of each triples $(\mathfrak{g},\mathfrak{h},\mathfrak{l}),$ $(\mathfrak{g},\mathfrak{l},\mathfrak{h})$ is described in Theorem \ref{twg} and given in Tables \ref{tab1} and \ref{tab2}.
\end{remark}
\begin{proof}
One easily verifies that in each case $\mathfrak{h}\cap \mathfrak{l}$ is of compact type and $\mathfrak{g}=\mathfrak{h}+\mathfrak{l}.$ Therefore we only need to show that there are no other possibilities. 

\textbf{Cases 1.-3.} contain all \textit{decomposable} triples, that is triples such that
$$\mathfrak{g}_{i} = \mathfrak{g}_{i}\cap \mathfrak{h} + \mathfrak{g}_{i}\cap \mathfrak{l}, \ \ i=1,2.$$
Thus decomposition of $\mathfrak{g}$ is obtained from decompositions of simple Lie algebras $\mathfrak{g}_{1}, \mathfrak{g}_{2}.$ But all possible proper semisimple decomposition of a simple Lie algebras (with $\mathfrak{h}\cap \mathfrak{l}$ of compact type) are given in Table \ref{tab2} (Theorem 4.1 and Table 2 in \cite{on1}). So the triples $(\mathfrak{g}_{i}, \mathfrak{g}_{i}\cap \mathfrak{h}, \mathfrak{g}_{i}\cap \mathfrak{l})$ are given in Table \ref{tab1} or are equal to $(\mathfrak{g}_{i}, \mathfrak{g}_{i}, \{ 0 \}).$

\textbf{Cases 4.-5.} one of subalgebras from the triple $(\mathfrak{g}, \mathfrak{h}, \mathfrak{l})$ contains a non-zero ideal of $\mathfrak{g} .$ Without loss of generality assume that $\mathfrak{g}_{1}\oplus \{ 0 \} \subset \mathfrak{h}.$ 
\begin{lemma}
The subalgebra $\pi_{2} (\mathfrak{h})\cap \pi_{2} (\mathfrak{l})$ is of compact type.
\end{lemma}
\begin{proof}
\textcolor[rgb]{0,0,1}{Choose a Cartan involution $\theta$ of $\mathfrak{g} .$ We may assume that $\theta (\mathfrak{h})=\mathfrak{h}$ and $\theta (\mathfrak{l})=\mathfrak{l} .$ We obtain a Cartan decomposition $\mathfrak{g}=(\mathfrak{k}_{1}\oplus \mathfrak{p}_{1})\oplus (\mathfrak{k}_{2} \oplus \mathfrak{p}_{2}).$ Assume that $X\in \pi_{2} (\mathfrak{h})\cap \pi_{2} (\mathfrak{l})$ does not lie in the maximal compact subalgebra of $\mathfrak{g} ,$ $\mathfrak{k}_{1}\oplus \mathfrak{k}_{2} .$ We see that $(0,X)\in \mathfrak{h}$ and there exists $(Y,0)\in \mathfrak{h}$ such that $(Y,X)\in \mathfrak{l} .$ But now $(Y,X)\in \mathfrak{h}\cap \mathfrak{l}$ and so $(X,Y) - (\theta (Y), \theta (X))\in \mathfrak{h}\cap \mathfrak{l}$ is a non-trivial element of $\mathfrak{p}_{1}\oplus \mathfrak{p}_{2}.$ A contradiction.}
\end{proof}
Since $\pi_{2}(\mathfrak{h})\cap \pi_{2}(\mathfrak{l})$ is of compact type and $\mathfrak{g}_{2}=\pi_{2}(\mathfrak{h})+\pi_{2}(\mathfrak{l})$ thus either $\pi_{2}(\mathfrak{h})$ is of compact type (and so equal to $\{ 0 \} $) or $(\mathfrak{g}_{2},\pi_{2}(\mathfrak{h}),\pi_{2}(\mathfrak{l}))$ is contained in Table \ref{tab1}. 

Notice that $\pi_{2}|_{\mathfrak{l}}$ is injective. Indeed, $\mathfrak{l}$ does not contain a non-trivial compact factor and
$$\operatorname{Ker} (\pi_{2}|_{\mathfrak{l}})= (\mathfrak{g}_{1}\oplus \{ 0 \} )\cap \mathfrak{l} \subset \mathfrak{h}\cap \mathfrak{l} .$$
Since the triple $(\mathfrak{g},\mathfrak{h},\mathfrak{l})$ is indecomposable and $\mathfrak{l}$ is simple so $\pi_{t}(\mathfrak{l})\cong \mathfrak{l},$ $t=1,2.$
 
\textbf{Case 6.} describes all other triples. In \cite{on1, on2} such triples are called primitive indecomposable and the classification of the triples is done using principal schemes. In our setting (where $\mathfrak{g}$ is a sum of two simple ideals), Theorem 4.3 in \cite{on1} can be stated as follows
\begin{proposition}
If the triple $(\mathfrak{g},\mathfrak{h},\mathfrak{l})$ is indecomposable and neither $\mathfrak{h}$ nor $\mathfrak{l}$ contains a non-zero ideal of $\mathfrak{g}$ then there is a triple $(\mathfrak{g}_{1},\mathfrak{g}',\mathfrak{g}'')$ from Table \ref{tab1} and a subalgebra $\mathfrak{g}_{2} \subset \mathfrak{g},$ $\mathfrak{g}''\subset \mathfrak{g}_{2}\subset \mathfrak{g}_{1},$ such that (after switching $\mathfrak{h}$ and $\mathfrak{l},$ if necessary) $\mathfrak{h}=\mathfrak{g}'\oplus \mathfrak{g}'',$ $\mathfrak{l}=\mathfrak{g}_{2}$ and
	$$\pi_{1}(\mathfrak{l})\neq \{ 0 \}, \ \pi_{2}(\mathfrak{l}) = \mathfrak{g}_{2} .$$
\end{proposition}

\noindent
This completes the proof of Theorem \ref{twg}.
\end{proof}

\textit{Proof of Proposition \ref{prop}.} Consider the space $SO(4,4)\times SO(2,4)/\Delta (SO(2,4))$ (the proof for $SO(3,4)\times SO(2,4)/\Delta (SO(2,4))$ is similar). Take $L'=SO(1,4),$ $L=Spin(3,4)\times L'.$ Let $\Gamma_{L'}$ be a torsion-free co-compact lattice of $L'$ and $\tilde{\Gamma}\subset Spin(3,4)$ a torsion-free co-compact lattice of $Spin(3,4).$ Since $(G\times H)/\textrm{diag}H\cong (G\times G)/\textrm{diag}(G)$ thus by Proposition \ref{propkako} there exists a neighborhood $U_{\epsilon}\subset \textrm{Hom}(\Gamma_{L'}, SO(2,4))$ of the natural inclusion such that for any $\rho\in U_{\epsilon},$ $\tilde{\Gamma}\times \rho (\Gamma_{L'})$ is a discrete subgroup that acts properly and co-compactly on $SO(4,4)\times SO(2,4)/\Delta (SO(2,4)).$ By Lemma 6.2 in \cite{kas} we can chose $\rho\in U_{\epsilon}$  so that $\rho (\Gamma_{L'})$ is Zariski dense in $SO(2,4).$  This way we obtain a discrete subgroup $\tilde{\Gamma}\times \rho (\Gamma_{L'})$ that acts properly and co-compactly on $SO(4,4)\times SO(2,4)/\Delta (SO(2,4)) $ with $\rho (\Gamma_{L'})$ Zariski dense in $SO(2,4).$ Assume that there exists a reductive subgroup $M\subset SO(4,4)\times SO(2,4)$ such that $M$ acts properly and co-compactly on $SO(4,4)\times SO(2,4)/\Delta (SO(2,4)) $ and $\tilde{\Gamma}\times \rho (\Gamma_{L'})\subset M.$ As $\rho (\Gamma_{L'})$ is Zariski dense in $SO(2,4)$ thus $Spin(3,4)\times SO(2,4)\subset M.$ 

\noindent
Let $\tilde{G}$ be an arbitrary Lie group with finitely many connected components and $\tilde{K}\subset\tilde{G}$ its maximal compact subgroup. Put $d(\tilde{G}):=\operatorname{dim}\tilde{G}-\operatorname{dim}\tilde{K}.$ We have
\begin{theorem}[\cite{ko-ma}, Theorem 4.7]
If $(G,H,L)$ corresponds to a standard compact Clifford-Klein form then $d(G)=d(H)+d(L).$
\end{theorem}

\noindent
But $d(SO(4,4)\times SO(2,4))=24,$ $d(M)\geq d(Spin(3,4)\times SO(2,4))=20,$ $d(SO(2,4))=8.$ A contradiction. This way we obtain a non-standard compact Clifford-Klein form of $SO(4,4)\times SO(2,4)/\Delta (SO(2,4)).$

\section{Tables}

Table \ref{tab1} contains all possible proper semisimple decompositions of absolutely simple real Lie algebras $\mathfrak{g},$ such that $\mathfrak{g}=\mathfrak{g}'+\mathfrak{g}''$ and $\mathfrak{g}'\cap\mathfrak{g}''$ is of compact type (see \cite{on1}, Theorem 4.1, Table 2 and \cite{koinv}, Section 5).

\begin{center}
 \begin{table}[h]
 \centering
 {\footnotesize
 \begin{tabular}{| c | c | c | c | c |}
   \hline
   \multicolumn{5}{|c|}{ \textbf{\textit{Decompositions of absolutely simple $\mathfrak{g}$}}} \\
   \hline                        
   $\mathfrak{g}$ & $\mathfrak{g}'$ & $\textcolor[rgb]{0,0,1}{\iota'}$ &  $\mathfrak{g}''$ & $\textcolor[rgb]{0,0,1}{\iota''}$ \\
   \hline
   $\mathfrak{su}(2,2n)$    & $\mathfrak{sp}(1,n)$ & $\xymatrix@1@R=2pt@!C=3pt{ 1 & && & & \\ {\circ}   \ar@{-}[r]& {\circ}   \ar@{-}[r]& \ar@{.}[r]&\ar@{-}[r] &{\circ} \ar@{<=}[r]&{\circ}}$ &   $\mathfrak{su}(1,2n) $ & $\xymatrix@1@R=2pt@!C=3pt{ 1 & && &  \\ {\circ}   \ar@{-}[r]& {\circ}   \ar@{-}[r]& \ar@{.}[r]&\ar@{-}[r] &{\circ}} $ $+ N$ \\
   \hline
   $\mathfrak{so}(2,2n)$     & $\mathfrak{su}(1,n)$ &  $\xymatrix@1@R=2pt@!C=3pt{ 1 & && & 1 \\ {\circ}   \ar@{-}[r]& {\circ}   \ar@{-}[r]& \ar@{.}[r]&\ar@{-}[r] &{\circ}} $   & $\mathfrak{so}(1,2n) $ & $\xymatrix@1@R=2pt@!C=3pt{ 1 & && & & \\ {\circ}   \ar@{-}[r]& {\circ}   \ar@{-}[r]& \ar@{.}[r]&\ar@{-}[r] &{\circ} \ar@{=>}[r]&{\circ}} $ $ + N $ \\
	& &  $+$ $\xymatrix@1@R=2pt@!C=3pt{  & && & 1 \\ {\circ}   \ar@{-}[r]& {\circ}   \ar@{-}[r]& \ar@{.}[r]&\ar@{-}[r] &{\circ}}$ & & \\
   \hline
   $\mathfrak{so}(4,4n)$   & $ \mathfrak{so}(3,4n)$ & $\xymatrix@1@R=2pt@!C=3pt{ 1 & && & & \\ {\circ}   \ar@{-}[r]& {\circ}   \ar@{-}[r]& \ar@{.}[r]&\ar@{-}[r] &{\circ} \ar@{=>}[r]&{\circ}} $ $+ N $  & $ \mathfrak{sp}(1,n)$ & $\xymatrix@1@R=2pt@!C=3pt{ 1 & && & & \\ {\circ}   \ar@{-}[r]& {\circ}   \ar@{-}[r]& \ar@{.}[r]&\ar@{-}[r] &{\circ} \ar@{<=}[r]&{\circ}} $  \\
	& & & & $+$ $\xymatrix@1@R=2pt@!C=3pt{ 1 & && & & \\ {\circ}   \ar@{-}[r]& {\circ}   \ar@{-}[r]& \ar@{.}[r]&\ar@{-}[r] &{\circ} \ar@{<=}[r]&{\circ}}$ \\
   \hline
   $\mathfrak{so}(4,4)$    & $ \mathfrak{so}(1,4)$ & $\xymatrix@1@R=2pt@!C=3pt{ 1  & \\ {\circ}   \ar@{=>}[r]& {\circ}} $ $ + 3N $  & $\mathfrak{so}(3,4)$ &  $\xymatrix@1@R=2pt@!C=3pt{  & & 1 \\ {\circ}   \ar@{-}[r]& {\circ}   \ar@{=>}[r]& {\circ}} $ \\
   \hline
   $\mathfrak{so}(3,4)$   & $ \mathfrak{so}(1,4)$ & $\xymatrix@1@R=2pt@!C=3pt{ 1  & \\ {\circ}   \ar@{=>}[r]& {\circ}} $ $ + 2N $  & $\mathfrak{g}_{2(2)}$ & $\xymatrix@1@R=2pt@!C=3pt{ &1&  \\ {\circ}  \ar@3{->}[r]&{\circ}   }$ \\
   \hline
	 $\mathfrak{so}(8,8)$   & $\mathfrak{so}(7,8)$ & $\xymatrix@1@R=2pt@!C=3pt{ 1 & && & & \\ {\circ}   \ar@{-}[r]& {\circ}   \ar@{-}[r]& \ar@{.}[r]&\ar@{-}[r] &{\circ} \ar@{=>}[r]&{\circ}} $ $+ N $  & $\mathfrak{so}(1,8)$ & $\xymatrix@1@R=2pt@!C=3pt{ & & & 1 \\ {\circ}   \ar@{-}[r]& {\circ}   \ar@{-}[r]& {\circ}   \ar@{=>}[r]& {\circ}} $ \\
   \hline 
 \end{tabular}
 }
 \caption{
 Proper semisimple decompositions of absolutely simple real Lie algebras inducing standard compact Clifford-Klein forms. In the Table $i',$ $i''$ denote the embeddings of $\mathfrak{g}',$ $\mathfrak{g}''$ into $\mathfrak{g},$ respectively (in the corresponding columns linear representations of the complexified Lie algebras $\mathfrak{g}'^{c},$ $\mathfrak{g}''^{c} $  realizing these embeddings are given and $N$ denotes the null representation). }
 \label{tab1}
 \end{table}
\end{center}

\noindent
Table \ref{tab2} contains triples described in Case 6. of Theorem \ref{twg}.

\begin{center}
 \begin{table}[h]
 \centering
 {\footnotesize
 \begin{tabular}{| c | c | c |}
   \hline
   \multicolumn{3}{|c|}{ \textbf{\textit{Case 6. triples}}} \\
   \hline                        
   $\mathfrak{g}$ & \textcolor[rgb]{0,0,1}{$\mathfrak{h}$} & \textcolor[rgb]{0,0,1}{$\mathfrak{l}$} \\
   \hline
   $\mathfrak{su}(2,2n)\oplus \mathfrak{su}(2,2n)$  & $\mathfrak{sp}(1,n) \oplus \mathfrak{su}(1,2n)$ & $\mathfrak{su}(2,2n)$ \\
   \hline
   $\mathfrak{so}(2,2n)\oplus \mathfrak{so}(2,2n)$  & $\mathfrak{su}(1,n) \oplus \mathfrak{so}(1,2n)$ & $\mathfrak{so}(2,2n)$ \\
   \hline
   $\mathfrak{so}(4,4n)\oplus \mathfrak{so}(4,4n)$  & $\mathfrak{so}(3,4n)\oplus \mathfrak{sp}(1,n)$ & $\mathfrak{so}(4,4n)$ \\
   \hline
   $\mathfrak{so}(4,4)\oplus \mathfrak{so}(4,4)$  & $\mathfrak{so}(3,4)\oplus \mathfrak{so}(1,4)$ & $\mathfrak{so}(4,4)$ \\
   \hline
	 $\mathfrak{so}(4,4)\oplus \mathfrak{so}(3,4)$  & $\mathfrak{so}(3,4)\oplus \mathfrak{so}(1,4)$ & $\mathfrak{so}(3,4)$ \\
   \hline
	 $\mathfrak{so}(4,4)\oplus \mathfrak{so}(2,4)$  & $\mathfrak{so}(3,4)\oplus \mathfrak{so}(1,4)$ & $\mathfrak{so}(2,4)$ \\
   \hline
   $\mathfrak{so}(3,4)\oplus \mathfrak{so}(3,4)$  & $\mathfrak{g}_{2(2)}\oplus \mathfrak{so}(1,4)$ & $\mathfrak{so}(3,4)$ \\
   \hline
	 $\mathfrak{so}(3,4)\oplus \mathfrak{so}(2,4)$  & $\mathfrak{g}_{2(2)}\oplus \mathfrak{so}(1,4)$ & $\mathfrak{so}(2,4)$ \\
   \hline
	 $\mathfrak{so}(8,8)\oplus \mathfrak{so}(8,8)$ & $\mathfrak{so}(7,8)\oplus \mathfrak{so}(1,8)$ & $ \mathfrak{so}(8,8) $ \\
   \hline 
 \end{tabular}
 }
\caption{
 Triples corresponding to Case 6. in Theorem \ref{twg}.}
 \label{tab2}
 \end{table}
\end{center}

\end{document}